# On three consecutive primes

Tsutomu Hashimoto

September, 2009

## Abstract

In this paper, we prove certain theorems about three consecutive primes.

## Introduction

The following result is due to Ishikawa [5].

**Theorem 1.**

If $p_n$ is the $n$th prime, then $p_{n+2} \leq p_{n+1} + p_n$ for all $n$.

We prove stronger results.

**Theorem 2.**

If $n \geq 4$, the inequality $p_{n+2}^2 < p_{n+1}^2 + p_n^2$ holds.

**Theorem 3.**

If $n \geq 9$, the inequality $p_{n+2}^3 < p_{n+1}^3 + p_n^3$ holds.

**Theorem 4.**

Let $k$ and $\varepsilon$ be any real numbers with $\varepsilon > 0$. Then
$$(2 - \varepsilon)p_{n+2}^k < p_{n+1}^k + p_n^k < (2 + \varepsilon)p_{n+2}^k$$
for all $n \geq$ some $n(k, \varepsilon)$.

## 1. Proof of Theorem 1

According to "A Theorem of Sylvester and Schur" [2] (see also [1]), there is a number containing a prime divisor $> p_n$ in the sequence $p_{n+1} + 1, p_{n+1} + 2, \ldots, p_{n+1} + p_n$. It follows, since there's not a multiple of $p_{n+1}$ in this sequence, that $p_{n+2} \leq p_{n+1} + p_n$. ∎

## 2. Proof of Theorem 2

According to Nagura [4] (see also [3]), $p_{n+1} - p_n - 1 < (1/5)p_n$ for $n > 9$.
Hence by a calculation
$$p_{n+1}^2 + p_n^2 - p_{n+2}^2 > (1/5^4)(229p_n^2 - 2460p_n - 2400) > 0 \; for \; n > 9.$$
By actual verification, we find that it is true for smaller values. ∎

## 3. Proof of Theorem 3

According to Rohrbach & Weis [6] (see also [3]), $p_{n+1} - p_n - 1 < (1/13)p_n$ for $n > 118$. Hence by a calculation

$$p_{n+1}^3 + p_n^3 - p_{n+2}^3 > (1/13^6)(3325841 p_n^3 - 23658180 p_n^2 - 56847882 p_n - 38416742) > 0 \text{ for } n > 118.$$

By actual verification, we find that it is true for smaller values. ∎

## 4. Proof of Theorem 4

Suppose that the notation $a_n^k$ stands for the sequence $((p_{n+1}/p_{n+2})^k + (p_n/p_{n+2})^k)$. It follows, since $\lim_{n\to\infty} a_n^k = \lim_{n\to\infty}((p_{n+1}/p_{n+2})^k + (p_n/p_{n+2})^k) = (1^k + 1^k) = 2$, that $|a_n^k - 2| < \varepsilon$ for all $n \geq$ some $n(k, \varepsilon)$. Hence

$$(2 - \varepsilon)p_{n+2}^k < p_{n+1}^k + p_n^k < (2 + \varepsilon)p_{n+2}^k,$$

because $|(p_{n+1}/p_{n+2})^k + (p_n/p_{n+2})^k - 2| < \varepsilon$. ∎

### Acknowledgement

I thank Jonathan Sondow for pointing out Ishikawa's paper and for suggesting Theorems 3 and 4.

Shiga 520-2412 JAPAN

t-hashimoto@aquablue.ne.jp